\documentclass[11pt]{article}
\def\NamedDraft{}
\usepackage[T1]{fontenc}
\usepackage{lmodern}
\usepackage[margin=1in]{geometry}
\usepackage{amsmath,amssymb,amsthm,mathtools}
\usepackage{graphicx,tikz,booktabs,enumitem,setspace,microtype,xcolor}
\usepackage[colorlinks=true,linkcolor=blue!40!black,
  citecolor=blue!40!black,urlcolor=blue!40!black]{hyperref}
\usetikzlibrary{arrows.meta}
\newtheorem{theorem}{Theorem}[section]
\newtheorem{proposition}[theorem]{Proposition}
\newtheorem{lemma}[theorem]{Lemma}
\newtheorem{corollary}[theorem]{Corollary}
\theoremstyle{remark}
\newtheorem{remark}[theorem]{Remark}
\newcommand{\E}{\mathbb E}
\newcommand{\T}{\mathcal T}

\newcommand{\depth}{\operatorname{depth}}

\definecolor{editcolor}{HTML}{C00030}

\setlist[enumerate]{itemsep=2pt,topsep=4pt}
\definecolor{SourceBlue}{HTML}{006BA4}
\definecolor{SupplyRed}{HTML}{C23B35}
\definecolor{DemandCyan}{HTML}{008B9B}
\definecolor{CellOne}{HTML}{006BA4}
\definecolor{CellTwo}{HTML}{5F9E6E}
\definecolor{CellThree}{HTML}{E69F00}
\definecolor{CellFour}{HTML}{986BB6}
\definecolor{CellFive}{HTML}{4BA3A8}
\definecolor{CellSix}{HTML}{D47B64}
\definecolor{CellSeven}{HTML}{8C8C8C}
\definecolor{CellEight}{HTML}{CC79A7}

\title{Near-Clairvoyant Dynamic Two-Sided Matching on the Line}
\ifdefined\NamedDraft
  \author{Yash Kanoria\\
    \small Columbia Business School\\[-2pt]
    \small \href{mailto:ykanoria@columbia.edu}{ykanoria@columbia.edu}}
\else
  \author{}
\fi
\date{}
\hypersetup{pdftitle={Near-clairvoyant dynamic two-sided matching on the line}}

\begin{document}
\maketitle

\begin{abstract}
Suppose $m$ supply units are available on the unit interval. At each time unit, a demand unit arrives at a 
uniformly random location on the interval, and must be matched to one of \(m\) available supply units, 
incurring a cost equal to the distance between the demand and supply locations. 
A new supply unit then arrives at an independent uniformly random location on the interval,
bringing the total number of supply units back to $m$.
The goal is to minimize the long run expected average match distance. 
We compare achievable cost against that achievable under perfect foresight; the latter is known to scale as $\Omega(\log m/m)$. %
Previously, an \(O(\sqrt{\log m})\)-factor online approximation to the clairvoyant cost was known. 
We construct an efficient online policy for this problem which approximates the clairvoyant cost to within an $O(1)$ factor, %
by managing stochastic fluctuations in supply inventory at all length scales, while enabling cancellation across scales. %
\end{abstract}
\noindent\textbf{Keywords:} dynamic two-sided matching; stochastic control;
spatial matching; Lyapunov drift method; online algorithms; hard \(k\)-taxi problem.

\section{Introduction}\label{sec:intro}

An operator matching demand and supply over time must decide which
available supply to use without knowing where future demand and supply will arise.
A nearby match saves distance now, but can leave too little supply where
it will next be needed. How costly is this uncertainty?

We investigate this question by studying the following model: 
Initially there are \(m\ge2\) supply units at arbitrary given
locations in \([0,1]\). At each period, the following events
occur in sequence:
\begin{itemize}
  \item A demand unit arrives at a uniformly random location in
  \([0,1]\).
  \item The operator chooses one of the available supply units to
  match to the demand. The matched pair leaves, and the period cost
  is the distance between their locations.
  \item A supply unit arrives at a uniformly random location in
  \([0,1]\), restoring the number of available supply units to \(m\).
\end{itemize}
All demand and replenishment locations are mutually independent.
An online policy chooses each match using the current
demand, available supplies, and past observations, before seeing the
new supply location.

The goal is to minimize the expected average cost per match over a
horizon of \(T\) periods: the expected total matching distance
divided by \(T\). We ask how this cost scales with \(m\), and compare
it with the cost achievable under perfect foresight. A clairvoyant
policy knows all arrival locations in advance but must still match
each demand immediately to an available supply unit.

This model is the one-dimensional case of the ``fully dynamic'' model of spatial matching 
introduced by Kanoria \cite{Kanoria} where demand and supply are located uniformly at random
in $[0,1]^d$. 
\cite{Kanoria} established
order-optimal online policies for minimizing expected average match distance for $d\ge2$, 
which achieve the same scaling as the clairvoyant lower bound 
(and best-case nearest neighbor distance) of $\Theta(1/m^{1/d})$ in each case,
thus showing that uncertainty about the future only costs a constant factor for $d\geq 2$.
The scaling behavior of the ``hard'' case of $d=1$ remained open: \cite{Kanoria} showed that for $d=1$, even with perfect foresight, 
the best achievable cost is $\Omega(\log m/m)$, exceeding the best-case nearest neighbor distance scaling of $\Theta(1/m)$. 
The best known online policy for $d=1$ is due to Coffman, Gilbert, and Shor \cite{CGS},\footnote{We 
thank Alireza Amanihamedani  (and later, Yilun Chen) for
pointing out this connection. \cite{CGS} concerns
directed distance around a circle. Their
stationary upper bound therefore carries over to Euclidean 
matching cost on the unit interval, as pointed out by Amanihamedani: Under rotation-invariant stationarity,
a uniform-rotation argument bounds the rule's expected distance on the
interval by twice its expected directed circular distance.
Averaging over the demand location within each circular gap also shows
that the expected interval distance is at least two thirds of the
expected directed circular distance.} 
and achieves a cost of $\Theta((\log m)^{3/2}/m)$, and thus approximates the clairvoyant cost to within 
a $\Theta(\sqrt{\log m})$ factor, leaving open the possibility that uncertainty costs up to a $\sqrt{\log m}$ factor. Their policy is counterclockwise-only greedy, which does not enjoy cancellation across length scales, motivating us to search for better policies. 

In this paper we close the scaling behavior of the one-dimensional case of the model. 
We give a policy attaining expected distance and root-mean-square distance both of order \(\Theta(\log m/m)\) for sufficiently long horizons, matching
the clairvoyant lower bound.
Uncertainty thus costs at most a constant factor even in the one-dimensional case.

\begin{theorem}[Optimal matching distance]\label{thm:main}
There are universal constants \(c>0\) and $C< \infty$ such that, for every
\(m\ge2\), there is a deterministic online matching policy with
expected average matching distance at most \(C\log m/m\) for every
horizon \(T\ge m^2\). For every
horizon \(T\ge m^2 \), every policy has
expected average matching distance at least \(c\log m/m\), even with perfect
foresight.

For the proposed policy, the same upper bound holds for the expected
root-mean-square matching distance for every horizon $T \geq m^3$.
All bounds hold for arbitrary initial supply locations.
Expectations are over subsequent arrivals.
The policy uses \(O(\log m)\) real-arithmetic operations
per period and \(O(m)\) memory.
\end{theorem}

The lower bound is due to Kanoria \cite[Theorem~3 and its proof]{Kanoria};
its proof permits perfect foresight. Our contribution is the
matching upper bound and the stronger second-moment guarantee.
The policy is defined in Section~\ref{sec:policy}.

Our result also applies to the hard \(k\)-taxi problem, which asks how
to dispatch \(k\) taxis to serve ride requests while minimizing empty
travel (with travel times ignored). Section~\ref{sec:taxi} states the
resulting near-clairvoyant guarantee for independent uniform pickups
and destinations.

\paragraph{The connection with static matching in two dimensions.}
\cite{Kanoria} shows the $\Omega(\log m/m)$ scaling of cost under perfect foresight by establishing a 
close connection between the problem with perfect foresight and the classical two-dimensional 
static matching problem of Ajtai, Koml\'os, and Tusn\'ady \cite{AKT}. The idea is that time in the dynamic model acts as an additional spatial dimension. This connection remains pivotal in 
the design of a near optimal dynamic matching policy for our setting, 
hence we provide a brief discussion of static matching in two dimensions before diving into our policy and analysis. 

Consider the  problem of matching $N$ supply units to $N$ demand units, both located i.i.d. uniformly on $[0,1]^2$, 
with cost being the expected average distance between matched pairs \cite{AKT}.  
Given a demand unit, the distance to its nearest neighbor scales as $1/\sqrt{N}$, 
and this is obviously a lower bound on the achievable cost. 
It turns out that every length scale $L \in [\Theta(1/\sqrt{N}), \Theta(1)]$ 
drives a lower bound with the same scaling: Consider an interior subsquare of side length $L$, 
and note that its area is $L^2$ and its perimeter is of order $L$. 
It has $NL^2 \pm \Theta(\sqrt{NL^2})$ demand units where the second term describes the stochastic uncertainty, 
and an i.i.d. number of supply units. Suppose, e.g., that realized demand in the subsquare exceeds supply by order $\sqrt{NL^2}= \sqrt{N} L$.
Then, in order to find enough supply to match all demand in the subsquare, one needs to include supply in a 
strip of width $\Delta_L \sim (\sqrt{N} L)/(NL) = 1/\sqrt{N}$ outside its perimeter, and one can show that, in fact, 
this constitutes a lower bound on the expected cost. 

For algorithmic purposes, it is natural to consider successive length scales which are a factor of 
two apart from each other, leading to order $\log N$ different length scales, each forcing a minimum expenditure of $1/\sqrt{N}$. 
If one is not careful about algorithm design, 
the costs at different length scales add up, leading to overall cost which scales as $\log N/\sqrt{N}$. This is 
a factor $\sqrt{\log N}$ worse than the optimal cost, which Ajtai et al.~\cite{AKT} show scales as $\sqrt{\log N/N}$. 
For instance, this $\sqrt{\log N}$ factor loss is known to occur under greedy-like algorithms \cite{Kanoria}. 
The optimal scaling can be achieved by arranging cancellation between displacements at different scales \cite{AKT}, which is challenging 
in a dynamic setting.
Returning to the one-dimensional dynamic model studied in this paper, %
the best known online algorithm \cite{CGS} is greedy like, 
and incurs an analogous $\Theta(\sqrt{\log m})$-factor loss relative to the clairvoyant lower bound. 
Similar to two dimensional static matching, efficient matching in this setting requires managing
stochastic fluctuations in supply inventory at all length scales, while enabling cancellation across scales
to eliminate the extra $\sqrt{\log m}$ factor.

\paragraph{Our policy and analysis.}
We now outline the main ideas of our policy and analysis. 

We organize supply into a dyadic tree of intervals, with levels of the tree corresponding to different length scales ranging from
$\Theta(1/m)$ to $\Theta(1)$, 
and nodes at each level corresponding to 
subintervals with the corresponding length. 
The state of the system is the supply counts at leaf nodes of the tree 
(and corresponding supply counts at each interior node of the tree).  
The policy is defined by a matching distribution over leaves of this tree 
(and corresponding distributions at each level of the tree); 
demand locations are matched to leaves via quantile matching.  
We use a logarithmic inventory potential capturing the evenness of the supply distribution, 
which sums over nodes at all levels of the tree with level weights chosen to ensure that levels (length scales) 
contribute about equally in the analysis (recall that the barrier to matching at all length scales is similar). 
The policy balances incurred match cost against balancing the supply distribution by appropriately biasing the matching distribution 
toward the more populous child at each node in the dyadic tree. Quantitatively, it uses harmonic feedback from the inventory state, corresponding to the discrete derivative of the potential function; an empty child receives no matching mass.

Tree symmetry leads to uncorrelated displacements at different levels of the tree (after an initial transient), 
facilitating cancellation across length scales. 
The analysis uses the Lyapunov drift method \cite{EryilmazSrikant} to compare the balancing effect of the policy 
with stochastic fluctuations. A Bellman correction accounts for how decisions at coarser 
scales affect inventory at finer scales. Averaging the resulting drift bound in steady state 
yields the sharp matching-distance guarantee, which is also shown to hold for long enough finite horizons. 

We present our policy and formal analysis in Sections~\ref{sec:policy}-\ref{sec:transport}, and return to motivate and discuss our choices in more detail in Section~\ref{sec:heuristic}.

\subsection{Related work}

\paragraph{Stochastic packing.}

A close conceptual precursor to the present paper, \cite{Shor1991} studies an online bin packing problem whose clairvoyant version relates to the two dimensional stochastic rightward matching problem, a variant of two dimensional static random matching \cite{AKT}. The achievability result is established by reducing the problem to a problem of putting balls in buckets on a ring, where each ball can go into one of two consecutive buckets. The construction of the online policy for the latter problem is based on a binary tree whose levels correspond to different length scales, with assignment biases defined at each node on the tree.  The tree construction leads to uncorrelated biases along the path from the root to any given node, and the squared discrepancy between the distribution of balls in buckets and the uniform distribution is studied and bounded. Different levels of the tree contribute equally (in terms of scaling) in the analysis. 

In the present paper, we also make use of a binary tree-based construction and uncorrelated biases along any path, and balance contributions from different levels of the tree. There are also some key differences, driven by differences between the two problem settings: we need to manage matching distance (as opposed to the hard constraint of having to choose one of two allowed buckets for each ball) while satisfying the hard constraint that all subintervals have non-negative inventory (as opposed to managing the discrepancy between bucket occupancies). 
We use imbalance-dependent bias magnitudes, whereas Shor only uses the sign of the local imbalance to determine the direction of the restoring bias. We study a quadratic ``energy'' which quantifies the expected squared matching distance, whereas Shor controls the squared discrepancy in the distribution of supply units (bucket occupancies). The logarithmic potential function and Lyapunov drift analysis we develop do not appear to have a clear analog in Shor's work.

\paragraph{Online spatial matching.}
Holden, Peres, and Zhai \cite{HPZ} consider $n$ supply units on the two 
dimensional surface of a three dimensional sphere, 
and introduce a gravitational allocation rule motivated by Newtonian mechanics 
to divide the sphere into $n$ cells of equal area. They show 
that this rule leads to a match distance which is within a constant factor
of the optimal transport under i.i.d. uniformly distributed locations. The allocation
rule applies naturally to online matching. 

In addition to the ``fully dynamic'' model studied in this paper, Kanoria \cite{Kanoria} introduced 
a ``semi-dynamic'' model of spatial matching with i.i.d. uniform demand and supply locations, where all supply is present initially
and is consumed without replacement. He provides online policies showing that uncertainty about the future 
costs only a constant factor in this model for all dimensions\footnote{The result for $d=2$ follows from \cite{HPZ}.} $d \in \mathbb{N}$, 
allowing for both equal supply and demand and excess supply. 
In the same model on the line, Balkanski, Faenza, and P\'erivier
\cite{BFP} prove constant competitiveness of nearest-available Greedy, both at balance
and with a fixed linear excess of supply.
Gupta et al.~\cite{GGPW} study online matching
with fixed supplies and i.i.d. requests from a known distribution over
the supply locations, without replenishment. They give a randomized
9-competitive algorithm on line and tree metrics.

Fielbaum et al. \cite{FCC} study our model of
replenished spatial matching (and generalizations of it) on a circle. Their results
describe the dynamics of the gap between consecutive supply units under greedy matching 
and structural properties of optimal assignments. 
They do not formally study the scaling of cost with the number of supply units.
Chen and Qian~\cite{ChenQian} also study our model on the circle. They show that two-sided greedy, which matches each demand to
the nearest available supply in either direction, unexpectedly incurs an average
expected cost per match of $\Omega(m^{-0.65})$, which is much
larger than the $\Theta((\log m)^{3/2}/m)$ cost of matching to the
first available supply in the counterclockwise direction~\cite{CGS}.

For adversarial online matching on the line with $n$ fixed supplies
and no replenishment, Raghvendra \cite{Raghvendra2018} gives a
deterministic $O(\log n)$-competitive algorithm, while Peserico and
Scquizzato \cite{PS2023} prove an $\Omega(\sqrt{\log n})$ lower bound
even for randomized algorithms. Thus, our result showing an $O(1)$ 
competitive ratio under i.i.d. uniform arrivals establishes that a natural ``average'' case 
behaves differently from the worst case. 

Chen et al.~\cite{FBDM} study ``feature-based'' dynamic matching with a fixed supply pool and sequential demand arrivals at i.i.d. locations. The demand and supply distributions are permitted to be different. 
One objective they study is to minimize the sum of squared matching distances.
Their policy simulates future demand to guide each match, and they bound its
regret per match relative to the limiting hindsight optimum.
 Kumar~\cite{KumarDark} obtains the same regret rates for broad classes of instances using only observed arrivals, without prior knowledge of the demand distribution.

Talagrand's  classical work on \emph{static} spatial matching in dimensions $d \geq 3$ of i.i.d. uniformly located supply and demand showed that the expected average match distance scales like the nearest neighbor distance \cite{Talagrand1992}. (Repeating the boundary strip calculation described above for the two-dimensional case provides intuition: for $d \geq 3$, the barrier to matching decays polynomially to $0$ in the length scale $L$, and hence the nearest neighbor distance dominates matching cost.)

\paragraph{$k$-taxi problems.}
In the $k$-taxi problem, $k$ taxis serve requests arriving one at a
time, each specifying a pickup and a destination. The selected taxi
travels empty to the pickup and then carries the passenger to the
destination. The ``easy'' version counts both empty and occupied
travel. Dehghani et al.\ \cite{Dehghani} study this objective under
known request distributions, calling it the ``Uber problem'',
and give constant-factor approximations to the optimal online
policy on the line and circle.

The ``hard'' version 
counts only empty travel. Since occupied
travel contributes an unavoidable cost shared by all policies,
a guarantee for total travel need not give a comparable guarantee
for empty travel. Competitive algorithms for adversarial hard
$k$-taxi have been developed by \cite{CK},
Gupta et al.~\cite{GKP}, and Coester et al.~\cite{CHP}.

Coester and Umenberger \cite{CU} obtain polylogarithmic competitive
guarantees for hard $k$-taxi in normed spaces when pickup and
destination densities remain bounded conditional on past requests.
Their model allows dependence across requests and between the
two endpoints. For independent uniform pickups and destinations
on the line, our result gives a constant competitive ratio and
identifies the optimal expected empty travel per request as
$\Theta(\log k/k)$; see Section~\ref{sec:taxi}.

\paragraph{The Lyapunov drift method.}
The Lyapunov drift method allows us to bound performance in stochastic systems by
relating the quantity of interest (such as queue length) to the \emph{drift} of a suitably chosen
potential (such as the sum of squared queue lengths), that is, its expected one-period change conditional on the
current state. It uses telescoping to show that the expected drift is zero (small on average) in steady state (over long horizons). Neely~\cite{Neely2010} uses drift inequalities to bound
time-average costs in stochastic networks. Eryilmaz and
Srikant~\cite{EryilmazSrikant} obtain bounds on mean queue lengths. In Section~\ref{sec:energy}, we use the
Lyapunov drift method to bound the ``transport energy'' of our matching
policy, leading to our bound on expected average match distance.

\section{The matching policy}\label{sec:policy}

Write \(S_t\) for the inventory after period \(t\), with \(S_0\)
the initial inventory. Let \(U_t\) be the demand, \(X_t\) the selected
supply, and \(r_t\triangleq|U_t-X_t|\) the matching distance in period \(t\).
When considering one period, write \(S,U,X,r\) without time subscripts.
All logarithms are natural unless a base is indicated. The policy in Theorem~\ref{thm:main} matches the first \(m\) demands
to the \(m\) initial supply units, retaining every
replenishment. At the start of period \(m+1\), the inventory thus consists
of \(m\) independent uniform supplies, and the policy switches to the
harmonic-feedback rule described below. Section~\ref{sec:finite} bounds
the cost of the first \(m\) matches.

\subsection{Matching probabilities on a dyadic tree}

The policy is defined for any dyadic number of cells \(n=2^{\ell_0}\) and any
regularization parameter \(a>0\). For the performance guarantee and its
proof, we choose\footnote{Our distance analysis can be adapted to any
dyadic \(n=\Omega(m/\log m)\).}
\begin{equation}\label{eq:parameters}
\ell_0\triangleq\lceil\log_2 m\rceil\, ,\qquad n\triangleq2^{\ell_0}\, ,\qquad a\triangleq48\ell_0\, .
\end{equation}
Thus \(m\le n<2m\).

Partition \([0,1]\) into \(n\) equal dyadic leaves,
and let \(\T\) be their complete binary tree, with root \(o\).
Let \(\T_\ell\) be the nodes at depth \(\ell\), with the root being at depth $0$. Use half-open cells,
except for the rightmost cell. A node \(v\) has interval \(I_v\),
length \(p_v\triangleq2^{-\depth(v)}\), and inventory \(N_v\).
Define its regularized inverse inventory by
\begin{equation}\label{eq:h-def}
h_v\triangleq\frac{p_v}{N_v+a/2}\, .
\end{equation}
This quantity is defined even at empty nodes.

Starting with deletion mass \(q_o\triangleq1\), assign masses recursively.
At an internal node, abbreviate its length and mass by \(p,q\), and
order its child inventories as \(x\ge y\). Their inverse inventories are
\[
h_+\triangleq\frac{p}{2x+a}\, ,\qquad h_-\triangleq\frac{p}{2y+a}\, .
\]
If both children are occupied, set
\begin{equation}\label{eq:policy}
b\triangleq a(h_--h_+)=\frac{2ap(x-y)}{(2x+a)(2y+a)}
\end{equation}
and give masses \((q+b)/2\) and \((q-b)/2\) to the larger and smaller
children, respectively. 
Thus, the bias $b$ incorporates harmonic feedback from the inventory state.  
If \(y=0\), set \(b\triangleq q\), giving all mass to the
larger child. This also assigns zero mass below an empty parent.
Restore the spatial labels after computing the split.

\subsection{Feasibility}

Define the normalized excess of arrivals over deletions by
\begin{equation}\label{eq:ht}
t_v\triangleq\frac{p_v-q_v}{a}\, .
\end{equation}
Here \(t_v\) is a node variable, distinct from the period index \(t\).

\begin{lemma}[Invariant and feasibility]\label{lem:invariants}
At every node,
\begin{equation}\label{eq:invariant-domain}
t_v\le\frac{h_v}{2}\, ,
\qquad\text{equivalently}\qquad q_v\ge N_vh_v\, .
\end{equation}
The policy assigns valid deletion probabilities, with positive mass to
every occupied node and zero mass to every empty node.
\end{lemma}
\begin{proof}
The two inequalities are equivalent because \(p=(N+a/2)h\).
At the root, \(t=0\). Suppose the invariant holds at a parent.
For two occupied children, \(h_\pm\ge h/2\) because \(N_\pm\le N\).
The definitions of \(t\) and \(b\) give
\[
t_+=\frac{t-h_-+h_+}{2}\le\frac{h_+}{2}\, ,\qquad
t_-=\frac{t+h_--h_+}{2}\le\frac{h_-}{2}\, .
\]
Here the inequalities use the induction hypothesis \(t\le h/2\le h_\pm\).
The resulting masses satisfy \(q_\pm\ge N_\pm h_\pm>0\) and sum to \(q\).
If \(y=0\), the empty child's invariant is equality. The other child
has inverse inventory \(h/2\) and mass \(q\ge Nh\ge N(h/2)\).
This argument includes an empty parent, where \(N=q=0\).
Induction proves feasibility and the invariant at every node in the tree.
\end{proof}

\subsection{Quantile implementation}
\label{subsec:quantile-implementation}

Order the supply cells \(I_1,\ldots,I_n\) from left to right. Partition
the demand interval \([0,1]\) into consecutive intervals
\(J_1,\ldots,J_n\) of lengths \(q_1,\ldots,q_n\). A demand in \(J_i\)
is matched to a supply in \(I_i\). %
Thus the policy amounts to quantile assignment up to the selection of
a supply unit within each supply cell.
Since demand is uniform, this selects cell \(I_i\) with probability
exactly \(q_i\). The choice of supply within the cell does not affect
the count dynamics.

For the distance analysis, associate each demand \(U\) with an auxiliary
location \(Y\) in its assigned supply cell. 
Let \(Q\) be the \emph{proxy supply distribution}, obtained by
spreading each cell's selection probability \(q_i\) uniformly
over \(I_i\), and let \(F_Q\) be its cumulative distribution function.
Set \(Y\triangleq F_Q^{-1}(U)\). This map sends \(J_i\) into \(I_i\). Thus the
actual matched supply \(X\) and the auxiliary location \(Y\) lie in the
same cell, giving
\[
  |X-Y|\le\frac1n\, .
\]

\begin{proposition}[Distance bound for the policy]\label{prop:quantile}
Conditional on the current inventory \(S\), the matching distance under
the policy above satisfies
\begin{equation}\label{eq:leaf-cost}
\E[r\mid S]
\le\frac1n+\int_0^1|F_Q(z)-z|\,dz\, .
\end{equation}
\end{proposition}
\begin{proof}
The triangle inequality and \(|X-Y|\le1/n\) give
\(r=|U-X|\le|U-Y|+1/n\). Since \(U\) is uniform and
\(Y=F_Q^{-1}(U)\), the one-dimensional quantile identity gives
\[
\E[|U-Y|\mid S]=\int_0^1|F_Q(z)-z|\,dz\, .
\]
Taking conditional expectations in the triangle inequality proves
\eqref{eq:leaf-cost}.
\end{proof}

\begin{figure}[htbp]
  \centering
  \resizebox{\linewidth}{!}{\begin{tikzpicture}[x=1cm,y=1cm,font=\small]
\path[use as bounding box] (.4,.8) rectangle (13.6,5.8);
\begin{scope}[yshift=-.25cm]
\node[anchor=west,font=\small] at (.6,5.75) {Supply counts};
\draw[black!25,line width=.4pt] (7.00000,5.77000)--(3.80000,5.28000);
\draw[black!25,line width=.4pt] (7.00000,5.77000)--(10.20000,5.28000);
\draw[black!25,line width=.4pt] (3.80000,5.28000)--(2.20000,4.79000);
\draw[black!25,line width=.4pt] (3.80000,5.28000)--(5.40000,4.79000);
\draw[black!25,line width=.4pt] (10.20000,5.28000)--(8.60000,4.79000);
\draw[black!25,line width=.4pt] (10.20000,5.28000)--(11.80000,4.79000);
\draw[black!25,line width=.4pt] (2.20000,4.79000)--(1.40000,4.30000);
\draw[black!25,line width=.4pt] (2.20000,4.79000)--(3.00000,4.30000);
\draw[black!25,line width=.4pt] (5.40000,4.79000)--(4.60000,4.30000);
\draw[black!25,line width=.4pt] (5.40000,4.79000)--(6.20000,4.30000);
\draw[black!25,line width=.4pt] (8.60000,4.79000)--(7.80000,4.30000);
\draw[black!25,line width=.4pt] (8.60000,4.79000)--(9.40000,4.30000);
\draw[black!25,line width=.4pt] (11.80000,4.79000)--(11.00000,4.30000);
\draw[black!25,line width=.4pt] (11.80000,4.79000)--(12.60000,4.30000);
\node[fill=white,text=black,inner sep=1.5pt,font=\small] at (7.00000,5.77000) {$ 32 $};
\node[fill=white,text=black,inner sep=1.5pt,font=\small] at (3.80000,5.28000) {$ 20 $};
\node[fill=white,text=black,inner sep=1.5pt,font=\small] at (10.20000,5.28000) {$ 12 $};
\node[fill=white,text=black,inner sep=1.5pt,font=\small] at (2.20000,4.79000) {$ 8 $};
\node[fill=white,text=black,inner sep=1.5pt,font=\small] at (5.40000,4.79000) {$ 12 $};
\node[fill=white,text=black,inner sep=1.5pt,font=\small] at (8.60000,4.79000) {$ 6 $};
\node[fill=white,text=black,inner sep=1.5pt,font=\small] at (11.80000,4.79000) {$ 6 $};
\node[fill=white,text=black,inner sep=1.5pt,font=\small] at (1.40000,4.30000) {$ 2 $};
\node[fill=white,text=black,inner sep=1.5pt,font=\small] at (3.00000,4.30000) {$ 6 $};
\node[fill=white,text=black,inner sep=1.5pt,font=\small] at (4.60000,4.30000) {$ 4 $};
\node[fill=CellFour!20,text=CellFour,inner sep=1.5pt,font=\small] at (6.20000,4.30000) {$ 8 $};
\node[fill=white,text=black,inner sep=1.5pt,font=\small] at (7.80000,4.30000) {$ 5 $};
\node[fill=white,text=black,inner sep=1.5pt,font=\small] at (9.40000,4.30000) {$ 1 $};
\node[fill=white,text=black,inner sep=1.5pt,font=\small] at (11.00000,4.30000) {$ 3 $};
\node[fill=white,text=black,inner sep=1.5pt,font=\small] at (12.60000,4.30000) {$ 3 $};
\fill[CellOne!14] (0.60000,3.40) rectangle (2.20000,4.06);
\draw[black!35,line width=.35pt] (0.60000,3.40) rectangle (2.20000,4.06);
\node[text=black!65,font=\small] at (1.40000,3.18) {A};
\fill[CellTwo!14] (2.20000,3.40) rectangle (3.80000,4.06);
\draw[black!35,line width=.35pt] (2.20000,3.40) rectangle (3.80000,4.06);
\node[text=black!65,font=\small] at (3.00000,3.18) {B};
\fill[CellThree!14] (3.80000,3.40) rectangle (5.40000,4.06);
\draw[black!35,line width=.35pt] (3.80000,3.40) rectangle (5.40000,4.06);
\node[text=black!65,font=\small] at (4.60000,3.18) {C};
\fill[CellFour!30] (5.40000,3.40) rectangle (7.00000,4.06);
\draw[black!35,line width=.35pt] (5.40000,3.40) rectangle (7.00000,4.06);
\node[text=black!65,font=\small] at (6.20000,3.18) {D};
\fill[CellFive!14] (7.00000,3.40) rectangle (8.60000,4.06);
\draw[black!35,line width=.35pt] (7.00000,3.40) rectangle (8.60000,4.06);
\node[text=black!65,font=\small] at (7.80000,3.18) {E};
\fill[CellSix!14] (8.60000,3.40) rectangle (10.20000,4.06);
\draw[black!35,line width=.35pt] (8.60000,3.40) rectangle (10.20000,4.06);
\node[text=black!65,font=\small] at (9.40000,3.18) {F};
\fill[CellSeven!14] (10.20000,3.40) rectangle (11.80000,4.06);
\draw[black!35,line width=.35pt] (10.20000,3.40) rectangle (11.80000,4.06);
\node[text=black!65,font=\small] at (11.00000,3.18) {G};
\fill[CellEight!14] (11.80000,3.40) rectangle (13.40000,4.06);
\draw[black!35,line width=.35pt] (11.80000,3.40) rectangle (13.40000,4.06);
\node[text=black!65,font=\small] at (12.60000,3.18) {H};
\fill[SupplyRed] (1.00000,3.74) circle (.043);
\fill[SupplyRed] (1.80000,3.74) circle (.043);
\fill[SupplyRed] (2.33333,3.74) circle (.043);
\fill[SupplyRed] (2.60000,3.74) circle (.043);
\fill[SupplyRed] (2.86667,3.74) circle (.043);
\fill[SupplyRed] (3.13333,3.74) circle (.043);
\fill[SupplyRed] (3.40000,3.74) circle (.043);
\fill[SupplyRed] (3.66667,3.74) circle (.043);
\fill[SupplyRed] (4.00000,3.74) circle (.043);
\fill[SupplyRed] (4.40000,3.74) circle (.043);
\fill[SupplyRed] (4.80000,3.74) circle (.043);
\fill[SupplyRed] (5.20000,3.74) circle (.043);
\fill[SupplyRed] (5.50000,3.74) circle (.043);
\fill[SupplyRed] (5.70000,3.74) circle (.043);
\fill[SupplyRed] (5.90000,3.74) circle (.043);
\fill[SupplyRed] (6.10000,3.74) circle (.043);
\fill[SupplyRed] (6.30000,3.74) circle (.043);
\fill[SupplyRed] (6.50000,3.74) circle (.043);
\fill[SupplyRed] (6.70000,3.74) circle (.043);
\fill[SupplyRed] (6.90000,3.74) circle (.043);
\fill[SupplyRed] (7.16000,3.74) circle (.043);
\fill[SupplyRed] (7.48000,3.74) circle (.043);
\fill[SupplyRed] (7.80000,3.74) circle (.043);
\fill[SupplyRed] (8.12000,3.74) circle (.043);
\fill[SupplyRed] (8.44000,3.74) circle (.043);
\fill[SupplyRed] (9.40000,3.74) circle (.043);
\fill[SupplyRed] (10.46667,3.74) circle (.043);
\fill[SupplyRed] (11.00000,3.74) circle (.043);
\fill[SupplyRed] (11.53333,3.74) circle (.043);
\fill[SupplyRed] (12.06667,3.74) circle (.043);
\fill[SupplyRed] (12.60000,3.74) circle (.043);
\fill[SupplyRed] (13.13333,3.74) circle (.043);
\node[anchor=west,font=\small] at (.6,2.52) {Demand mapping under harmonic feedback};
\fill[CellOne!14] (0.60000,1.46) rectangle (1.79169,2.00);
\draw[black!45,line width=.45pt] (0.60000,1.46) rectangle (1.79169,2.00);
\node at (1.19584,1.73) {A};
\fill[CellTwo!14] (1.79169,1.46) rectangle (3.78338,2.00);
\draw[black!45,line width=.45pt] (1.79169,1.46) rectangle (3.78338,2.00);
\node at (2.78753,1.73) {B};
\fill[CellThree!14] (3.78338,1.46) rectangle (5.34459,2.00);
\draw[black!45,line width=.45pt] (3.78338,1.46) rectangle (5.34459,2.00);
\node at (4.56398,1.73) {C};
\fill[CellFour!30] (5.34459,1.46) rectangle (7.33247,2.00);
\draw[black!45,line width=.45pt] (5.34459,1.46) rectangle (7.33247,2.00);
\node at (6.33853,1.73) {D};
\fill[CellFive!14] (7.33247,1.46) rectangle (9.45887,2.00);
\draw[black!45,line width=.45pt] (7.33247,1.46) rectangle (9.45887,2.00);
\node at (8.39567,1.73) {E};
\fill[CellSix!14] (9.45887,1.46) rectangle (10.36623,2.00);
\draw[black!45,line width=.45pt] (9.45887,1.46) rectangle (10.36623,2.00);
\node at (9.91255,1.73) {F};
\fill[CellSeven!14] (10.36623,1.46) rectangle (11.88312,2.00);
\draw[black!45,line width=.45pt] (10.36623,1.46) rectangle (11.88312,2.00);
\node at (11.12468,1.73) {G};
\fill[CellEight!14] (11.88312,1.46) rectangle (13.40000,2.00);
\draw[black!45,line width=.45pt] (11.88312,1.46) rectangle (13.40000,2.00);
\node at (12.64156,1.73) {H};
\node[font=\scriptsize,text=black!65] at (.6,1.24) {$0$};
\node[font=\scriptsize,text=black!65] at (13.4,1.24) {$1$};
\draw[DemandCyan,densely dashed,line width=.5pt] (7.12800,2.00)--(7.12800,2.25);
\draw[DemandCyan,densely dashed,line width=.5pt] (7.12800,2.77)--(7.12800,3.40);
\draw[DemandCyan,line width=1.2pt] (6.50000,3.74)--(7.12800,3.74);
\draw[CellFour,line width=1.1pt] (6.50000,3.74) circle (.105);
\fill[DemandCyan] (7.12800,3.82200)--(7.21000,3.74000)--(7.12800,3.65800)--(7.04600,3.74000)--cycle;
\fill[DemandCyan] (7.12800,1.81200)--(7.21000,1.73000)--(7.12800,1.64800)--(7.04600,1.73000)--cycle;
\node[anchor=west,text=DemandCyan,font=\scriptsize] at (7.26800,2.88) {$u=0.51$};
\end{scope}
\end{tikzpicture}}
  \caption{Illustration of harmonic feedback with \(m=32\), \(n=8\),
  and \(a=4\), using illustrative parameters rather than those in
  \eqref{eq:parameters}. The lower intervals have lengths equal to
  the cells' matching probabilities. The demand at \(u=0.51\) is
  assigned to cell D and matched to the indicated supply.}
  \label{fig:policy}
\end{figure}
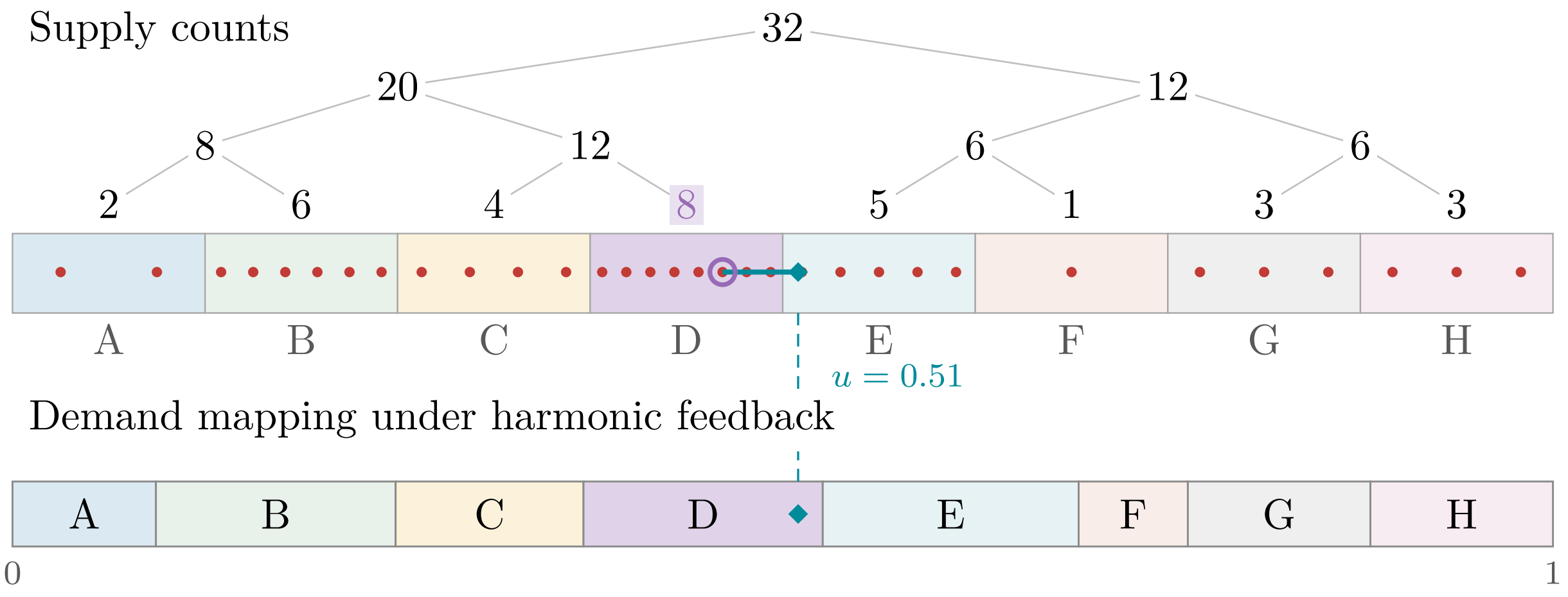

\paragraph{A numerical illustration.}
Figure~\ref{fig:policy} illustrates the general rule with \(m=32\),
\(n=8\), and \(a=4\). The cell counts are \((2,6,4,8,5,1,3,3)\).
At the root, \((x,y,p,q)=(20,12,1,1)\), so
\[
  b=\frac{2\cdot4\cdot1\cdot(20-12)}{(40+4)(24+4)}
   =\frac4{77}\, ,\qquad
  (q_{\mathrm{left}},q_{\mathrm{right}})
   =\left(\frac{81}{154},\frac{73}{154}\right)\, .
\]
The fuller half therefore has expected inventory change
\(1/2-81/154=-2/77\). Continuing down the tree assigns cell D,
whose physical interval is \([3/8,1/2]\), the demand interval
\([685/1848,81/154]\). Thus a demand at \(0.51\) uses a supply in D.
These illustrative parameters differ from those in \eqref{eq:parameters}:
for \(m=32\), the theorem uses \(n=32\) and \(a=240\).

\paragraph{Compute and memory requirements.}
For a given demand location \(U\), its assigned supply cell can be found
by following a single root-to-leaf path. Start at the root with \(s\triangleq0\)
and \(q\triangleq1\). At an internal node whose assigned demand interval is
\([s,s+q)\), use the policy rule to compute the left child's mass
\(q_{\mathrm L}\). Accordingly, follow the left child if
\(U<s+q_{\mathrm L}\), and the right child otherwise.
Store each cell's supplies in a linked list and select the first supply
in that list. After deletion and replenishment, update the counts along
the two affected paths. This requires \(O(\ell_0+1)\) operations per period
and \(O(m+n)\) memory.

\subsection{Symmetry of the count process}

\begin{lemma}[Count-chain symmetry]\label{lem:symmetry}
The finite chain of leaf counts is irreducible and aperiodic and hence has
a unique stationary law \(\pi\). This law is invariant under every
automorphism of the rooted binary tree. An invariant initial count law
remains invariant at every deterministic time.
\end{lemma}
\begin{proof}
Every occupied leaf has positive deletion probability, so supplies can be
moved between arbitrary leaves with positive probability. Replenishment in
the deleted leaf gives a self-loop. Thus the chain is irreducible and
aperiodic. Exchanging two sibling subtrees exchanges their counts and,
by \eqref{eq:policy}, their deletion masses. The transition kernel therefore
commutes with every tree automorphism. Uniqueness makes the stationary law
invariant, and the same kernel identity preserves any invariant initial law.
\end{proof}

\section{Bounding transport energy by the Lyapunov drift method}\label{sec:energy}

Write the signed bias \(b_v\triangleq q_{v_{\mathrm L}}-q_{v_{\mathrm R}}\) at node $v$ using the fixed
spatial child labels. Define the transport energy and the level energies by
\begin{equation}\label{eq:H-G}
G\triangleq\sum_{v\in\T\setminus\T_{\ell_0}}p_vb_v^2\, ,
\qquad H_\ell\triangleq\sum_{v\in\T_\ell}p_vh_v^2\, ,
\qquad H_0=\frac1{(m+a/2)^2}\, .
\end{equation}
Section~\ref{sec:transport} will use \(G\) to bound squared matching
distance.

To bound \(G\), we instantiate the Lyapunov drift method: we use a logarithmic inventory potential (Lyapunov function) \(\Phi\) that
rewards an even distribution of supply. We show that its conditional
expected increase in one period is at least a positive multiple of
\(G\) minus a state-independent term. Over \(T\) periods, the expected
increases telescope to \(\E[\Phi_T-\Phi_0]\). Since \(\Phi\) is
bounded, this gives an upper bound on the expected time average of
\(G\) from any initial inventory. 

From any initial inventory, Proposition~\ref{prop:energy} bounds
the expected time average of \(G\) over \(T\) periods by
\(20a^2/m^2+20a^2/(17T)\). In steady state, the bound is \(20a^2/m^2\).
The level energy \(H_\ell\) measures squared inverse inventory at depth
\(\ell\). Two local inequalities control \(G\) through a potential drift
and bound the increase of \(H_\ell\) down the tree.

\subsection{The harmonic potential}

Define the concave harmonic sum
\[
F_a(k)\triangleq\sum_{j=1}^k\frac1{j+a/2}\, ,
\qquad F_a(0)\triangleq0\, ,
\]
and define the inventory potential and its linearized drift
\begin{equation}\label{eq:Phi-D}
\Phi\triangleq\sum_{v\ne o}p_v^2F_a(N_v)\, ,
\qquad
D\triangleq\sum_{v\ne o}p_v^2\frac{p_v-q_v}{N_v+a/2}
=a\sum_{v\ne o}p_vt_vh_v\, .
\end{equation}
At each level of the tree, concavity makes its contribution to the potential \(\Phi\) largest
when counts are as equal as possible. The linearized drift \(D\) multiplies each
node's net inventory drift by its current harmonic increment \(1/(N_v+a/2)\). 
A net gain of one supply actually contributes the smaller
increment \(1/(N_v+a/2+1)\). We denote the resulting
{conditional expected curvature loss} by \(\mathcal R\).
The following lemma gives the exact drift identity.

\begin{lemma}[Exact harmonic drift]\label{lem:drift}
Conditional on the present inventory,
\begin{equation}\label{eq:drift}
\E[\Phi_{\mathrm{next}}-\Phi\mid S]=D-\mathcal R\, ,
\end{equation}
where
\begin{equation}\label{eq:curvature-error}
\mathcal R\triangleq
\sum_{v\ne o}
\frac{p_v^3(1-q_v)}{(N_v+a/2)(N_v+a/2+1)}
\le\sum_{\ell=1}^{\ell_0} H_\ell \, .
\end{equation}
\end{lemma}
\begin{proof}
At a node of count \(k\), abbreviate its spatial and deletion masses as
\(p,q\). The probabilities of an increase and a decrease are
\(p(1-q)\) and \(q(1-p)\). The conditional drift of its harmonic sum is
\[
\frac{p(1-q)}{k+a/2+1}-\frac{q(1-p)}{k+a/2}
=\frac{p-q}{k+a/2}
-\frac{p(1-q)}{(k+a/2)(k+a/2+1)}\, .
\]
This remains valid at \(k=0\), where \(q=0\).
Multiply by \(p_v^2\) and sum. The curvature contribution at a node is
at most \(p_vh_v^2\), proving the final inequality.
\end{proof}

\subsection{Two local inequalities}

In the local calculations, omit the parent index and order the child
counts as \(x\ge y\), with \(N\triangleq x+y\). The labels \(+,-\)
refer to the fuller and less full children, respectively, and \(b=|b_v|\) is the
magnitude of the spatially signed bias. 

A node inherits its total deletion mass \(q\) from its ancestors; its bias 
only determines how that mass is split between its children. The two child 
inventory drifts are \((p-q-b)/2\) and \((p-q+b)/2\). Their common component 
\((p-q)/2\) can be negative even when the children have equal counts. 
A quadratic correction accounts for this
inherited drift when local estimates are summed. A second inequality
controls the curvature loss across levels by the same squared bias.

Use the polynomial correction
\begin{equation}\label{eq:Bellman}
B(h,t)\triangleq th-\frac23t^2-\frac13h^2
=-\frac13(h-t)(h-2t)\, .
\end{equation}
Write \(B_v\triangleq B(h_v,t_v)\). The invariant \(t\le h/2\)
gives \(B\le0\).
For quantities attached to two children, \(\E_{\mathrm{ch}}\) denotes
their equal average, with spatial weights \(1/2,1/2\). Primes below
denote child quantities. %

\begin{lemma}[Local estimates]\label{lem:local}
At every internal node,
\begin{align}
\E_{\mathrm{ch}}[t'h'+B(h',t')]-B(h,t)
&\ge\frac{b^2}{30a^2}\, ,\label{eq:local}\\
\E_{\mathrm{ch}}(h')^2-h^2
&\le\frac{3b^2}{4a^2}\, .\label{eq:local-rate}
\end{align}
\end{lemma}

Before providing a proof, we indicate how the lemma will be used in the sequel. Inequality \eqref{eq:local} will be summed over internal nodes in the tree, causing the quadratic correction terms to telescope and leading to a lower bound on the linearized drift $D$. The inequality~\eqref{eq:local-rate} controls the curvature loss. Fix the current inventory. Multiplying the second inequality by \(p_v\)
and summing over nodes at depth \(\ell\) gives
\[
H_{\ell+1}-H_\ell
\le \frac{3}{4a^2}\sum_{v\in\T_\ell}p_vb_v^2\, .
\]
Summing over depths \(j=0,\dots,\ell-1\) yields \(H_\ell\le H_0+3G/(4a^2)\).
Lemma~\ref{lem:drift} bounds the one-period expected curvature loss by
\(\mathcal R\le\sum_{\ell=1}^{\ell_0} H_\ell\), so
\begin{equation}\label{eq:R-bound}
\mathcal R\le \ell_0H_0+\frac{3\ell_0}{4a^2}G\, .
\end{equation}

\begin{proof}[Proof of Lemma~\ref{lem:local}]
Set
\[
H\triangleq\frac{h_-+h_+}{2}\, ,\qquad w\triangleq\frac{h_--h_+}{2}\, ,\qquad
K\triangleq\frac{p}{N+a}\, .
\]
The harmonic mean of \(h_+\) and \(h_-\) is \(K\). Consequently,
\begin{equation}\label{eq:harmonic-identity}
w^2=H(H-K)\, ,\qquad H\ge K\, ,\qquad K\le h\le2K\, .
\end{equation}
These identities apply to empty children as well.

First suppose both children are occupied. Then \(b/a=2w\), and
\(h_\pm=H\mp w\), \(t_\pm=t/2\mp w\).
Write \(R\triangleq\E_{\mathrm{ch}}[t'h'+B(h',t')]-B(h,t)\).
Expanding \(B\) and averaging over the two children gives
\[
R=\frac{t^2}{2}+t(H-h)+\frac{h^2-H^2}{3}+w^2\, .
\]
Substituting \(w^2=H^2-HK\) from \eqref{eq:harmonic-identity}
and completing the square in \(t\) gives
\[
R-\frac{w^2}{6}
=\frac12(t+H-h)^2+\frac{(6h-5K)H-h^2}{6} \ge \frac{(h-K)(5K-h)}{6}\ge0 \, , 
\]
where the inequalities are obtained by discarding the nonnegative square and using \(H\ge K\) in the remaining term.
The final inequality uses \(K\le h\le2K\). Since \(b/a=2w\), we obtain
\[
R\ge\frac{w^2}{6}=\frac{b^2}{24a^2}\, .
\]

If a child is empty, the policy instead uses \(b=q\). Put
\[
P\triangleq\frac pa\, ,\qquad \rho\triangleq P-\frac h2\ge0\, ,\qquad
\eta\triangleq\frac h2-t\ge0\, ,\qquad
\beta\triangleq\frac ba=\frac qa=\rho+\eta\, .
\]
Here \(\eta\) is the slack in the parent invariant, and
\(\rho\ge0\) follows from \(h\le2p/a\).
The child variables are \((h_-,t_-)=(P,P/2)\) and
\((h_+,t_+)=(h/2,t-P/2)\). Direct expansion gives
\begin{equation}\label{eq:empty-square}
R=\frac{\beta^2+(3\eta-2\rho)^2}{30}+\frac{h^2}{8}+\frac{h\rho}{12}
\ge\frac{\beta^2}{30}\, .
\end{equation}
This includes an empty parent, for which \(\rho=\eta=\beta=0\), and proves
\eqref{eq:local}.

We now establish \eqref{eq:local-rate}. Using \eqref{eq:harmonic-identity} gives
\[
H^2-K^2=(H-K)(H+K)\le2H(H-K)=2w^2\, .
\]
Hence, using \(K\le h\),
\[
\E_{\mathrm{ch}}(h')^2=H^2+w^2
\le K^2+3w^2\le h^2+3w^2\, .
\]
For occupied children, \(2w=b/a\). For an empty child,
\(2w=\rho\le\rho+\eta=b/a\). This proves \eqref{eq:local-rate} in every case.
\end{proof}

\subsection{The energy bound}

The next proposition leverages the lower bound on $D$ and the upper bound on $\mathcal{R}$ implied by 
Lemma~\ref{lem:local} to obtain a lower bound on the expected change in the potential \eqref{eq:drift} 
in one time step. It then sums over these expected changes across time, and uses the fact that the potential 
is bounded to deduce an upper bound on the transport energy.

\begin{proposition}[Energy estimate]\label{prop:energy}
For arbitrary initial supply locations and every integer \(T\ge1\),
\begin{equation}\label{eq:finite-G}
\frac1T\sum_{t=1}^T\E G(t)
\le\frac{20a^2}{m^2}+\frac{20a^2}{17T}\, .
\end{equation}
Here \(G(t)\) is the energy before period \(t\).
Consequently, under the stationary distribution of the leaf counts,
\begin{equation}\label{eq:stationary-energy}
\E_\pi G\le\frac{20a^2}{m^2}\, .
\end{equation}
\end{proposition}
\begin{proof}
For a child quantity \(f\), the weighted average
\(p_v\E_{\mathrm{ch}}f\) equals the sum of \(p_wf_w\) over the
two children \(w\) of \(v\). Multiply Lemma~\ref{lem:local} inequality \eqref{eq:local} by \(p_v\)
and sum over internal nodes. Every nonboundary \(B_v\) then appears
once with each sign, leaving
\[
\frac Da+\sum_{v\in\T_{\ell_0}}p_vB_v-B_o\, .
\]
The leaf terms are nonpositive: Lemma~\ref{lem:invariants} gives
\(t_v\le h_v/2\), so the factorization \eqref{eq:Bellman} yields \(B_v\le0\).
Also, \(B_o=-H_0/3\), so
\begin{equation}\label{eq:aggregate}
D\ge\frac{G}{30a}-\frac a3H_0\, .
\end{equation}
\begin{samepage}
Combining \eqref{eq:aggregate} with the curvature bound \eqref{eq:R-bound}
in the drift identity \eqref{eq:drift}, and using \(a=48\ell_0\), gives
\begin{align*}
\E[\Phi_{\mathrm{next}}-\Phi\mid S]
&\ge\left(\frac1{30a}-\frac{3\ell_0}{4a^2}\right)G
-\left(\frac a3+\ell_0\right)H_0\\
&=\frac{17}{960a}G-\frac{17a}{48}H_0\, .
\end{align*}
\end{samepage}
Summing expected drifts over a finite horizon leaves the endpoint
difference \(\E[\Phi_T-\Phi_0]\), so
\begin{equation}\label{eq:master-energy}
\frac1T\sum_{t=1}^T\E G(t)
\le20a^2H_0+\frac{960a}{17T}\E[\Phi_T-\Phi_0]\, ,
\end{equation}
where \(\Phi_t\) is the potential after period \(t\).
Since \(x\mapsto1/(x+a/2)\) is decreasing,
\[
F_a(k)\le\int_0^k\frac{dx}{x+a/2}
=\log(1+2k/a)\, .
\]
Together with \(N_v\le m\) and
\(\sum_{v\ne o}p_v^2=\sum_{\ell=1}^{\ell_0}2^{-\ell}<1\),
this gives \(\Phi\le\log(1+2m/a)\).
It follows that 
\[
0\le\Phi\le\log(1+2m/a)\le\log m\le \ell_0=a/48\, ,
\]
using \(m\ge2\), and \(a\ge48\).
This proves \eqref{eq:finite-G}.
Taking \(T\to\infty\) in \eqref{eq:finite-G} proves
\eqref{eq:stationary-energy} under stationarity.
\end{proof}

\section{Bounding matching distance via cancellation across length scales}\label{sec:transport}

Recall the proxy supply distribution $Q$ from Section~\ref{subsec:quantile-implementation}. 
We now relate \(G\) to squared matching distance by expressing
\(F_Q(z)-z\) in terms of the split biases \(b_v\).
Since \(Q\) is uniform within each leaf, \(F_Q\) is linear there.
We recover \(F_Q\) by linear interpolation, starting with its
values at \(0\) and \(1\), then adding its values at the dyadic
midpoints one level at a time. The initial interpolation is
the function \(z\).

At the midpoint of an interval \(I_v\), the value of \(F_Q\)
differs from the average of its endpoint values by
\((q_{v_{\mathrm L}}-q_{v_{\mathrm R}})/2=b_v/2\).
Adding this midpoint therefore changes the interpolation by
\(b_v\tau_v\), where \(\tau_v\) is zero outside \(I_v\), linear
from zero at its left endpoint to \(1/2\) at its midpoint,
and linear back to zero at its right endpoint.
After all \(\ell_0\) levels, the interpolation equals \(F_Q\).
Summing these corrections gives
\begin{equation}\label{eq:Haar}
F_Q(z)-z=\sum_{v\in\T\setminus\T_{\ell_0}}b_v\tau_v(z)\, ,
\qquad
\int_0^1\tau_v(z)^2\,dz=\frac{p_v}{12}\, ,
\end{equation}
where the second identity follows by direct integration.
\if
We now relate \(G\) to squared matching distance. Integrating each
split's density deviation from the parent-cell average gives a tent.
For each internal node \(v\), let \(\tau_v\) be the tent supported on
\(I_v\), linear from zero at its left endpoint to height \(1/2\) at its
midpoint and back to zero at its right endpoint. With the signed spatial
bias \(b_v=q_{v_{\mathrm L}}-q_{v_{\mathrm R}}\), the deletion
distribution has the exact integrated Haar expansion
\begin{equation}\label{eq:Haar}
F_Q(z)-z=\sum_{v\in\T\setminus\T_{\ell_0}}b_v\tau_v(z)\, ,
\qquad
\int_0^1\tau_v(z)^2\,dz=\frac{p_v}{12}\, .
\end{equation}
\fi

\begin{lemma}[Cancellation and transport energy]\label{lem:transport-energy}
Under any tree-invariant count law,
\begin{equation}\label{eq:transport-energy}
\E\int_0^1(F_Q(z)-z)^2\,dz=\frac1{12}\E G\, .
\end{equation}
Consequently,
\begin{equation}\label{eq:transport-cost}
\sqrt{\E r^2}\le\frac1n+\sqrt{\frac{\E G}{12}}\, ,
\end{equation}
where $r\triangleq |U-X|$ is the match distance.
\end{lemma}
\begin{proof}
Tents belonging to disjoint nodes have disjoint interiors. If \(u\) is
an ancestor of \(v\ne u\), exchange the two child subtrees of \(v\).
This negates \(b_v\) and fixes \(b_u\). Tree invariance therefore gives
\(\E[b_ub_v]=0\). Thus the mixed terms for nested tents cancel in
expectation. Expanding the left-hand side \eqref{eq:transport-energy} using 
\eqref{eq:Haar} then establishes
\eqref{eq:transport-energy}.

For the second assertion, recall that \(U\) is the uniform demand location, 
\(Y\triangleq F_Q^{-1}(U)\) is the proxy supply location, 
and \(X\) is the matched supply location in the same supply 
cell as the proxy location (thus \(|X-Y|\le1/n\)).
Conditional on the inventory, \(Y\) has distribution \(Q\).
Write \(g(z)\triangleq F_Q(z)-z\). Continuity of \(F_Q\) gives
\(F_Q(Y)=U\) almost surely, hence \(Y-U=-g(Y)\). Using
\(F_Q'=1+g'\) gives
\begin{equation}
\label{eq:conditional-transport}
\E[(Y-U)^2\mid S]
=\int_0^1g(z)^2F_Q'(z)\,dz
=\int_0^1g(z)^2\,dz+\left[\frac{g(z)^3}{3}\right]_0^1
=\int_0^1g(z)^2\,dz\, .
\end{equation}
The boundary term is zero because \(g(0)=g(1)=0\).
As \(|X-Y|\le1/n\), the \(L^2\) triangle
inequality (the root mean square of a sum is at most the sum of the root mean squares) 
and \eqref{eq:transport-energy} prove
\eqref{eq:transport-cost}.
\end{proof}

\subsection{Proof of Theorem~\ref{thm:main}}\label{sec:finite}

First suppose the initial count law is tree invariant and run the
harmonic-feedback policy from the outset. Lemma~\ref{lem:symmetry}
preserves this invariance at every period. Apply
\eqref{eq:transport-cost} from Lemma~\ref{lem:transport-energy} at an independent uniformly selected time
among \(1,\ldots,T\), and use \eqref{eq:finite-G} from Proposition~\ref{prop:energy}, to obtain
\begin{equation}\label{eq:finite-cost}
 \biggl(\frac1T\sum_{t=1}^T\E r_t^2\biggr)^{1/2}
 \le \frac1n+a\sqrt{\frac5{3m^2}+\frac5{51T}}\, ,
 \qquad T\ge1\, .
\end{equation}

For an arbitrary initial inventory, first match one original supply
per period for \(m\) periods, retaining every
replenishment. Thereafter, use the harmonic-feedback policy.
Each of the first \(m\) distances is at most one, and the inventory
after those periods consists of \(m\) independent uniform supplies.
Its count law is tree invariant, so \eqref{eq:finite-cost} applies
to the remaining \(s\triangleq T-m\) periods. For \(T>m\), this gives
\[
 \biggl(\frac1T\sum_{t=m+1}^T\E r_t^2\biggr)^{1/2}
 \le \sqrt{\frac sT}
       \left(\frac1n+a\sqrt{\frac5{3m^2}+\frac5{51s}}\right)
 \le \frac1n+a\sqrt{\frac5{3m^2}+\frac5{51T}}\, .
\]
Cauchy--Schwarz also bounds the contribution of these periods to the
expected average distance by the last expression. Including the first
\(m\) periods therefore yields
\begin{align}
 \frac1T\sum_{t=1}^T\E r_t
 &\le \frac mT+\frac1n+
       a\sqrt{\frac5{3m^2}+\frac5{51T}}\, ,
 \label{eq:arbitrary-mean}\\
 \biggl(\frac1T\sum_{t=1}^T\E r_t^2\biggr)^{1/2}
 &\le \sqrt{\frac mT}+\frac1n+
       a\sqrt{\frac5{3m^2}+\frac5{51T}}\, ,
 \label{eq:arbitrary-rms}
\end{align}
where the second inequality uses \(\sqrt{x+y}\le\sqrt x+\sqrt y\).
For \(T\ge m^2\) in \eqref{eq:arbitrary-mean} and \(T\ge m^3\)
in \eqref{eq:arbitrary-rms}, these bounds are at most
\[
 \frac2m+\sqrt{\frac{30}{17}}\frac am
 <\frac{3a}{2m}=O(\log m/m)\, ,
\]
using \(n\ge m\) and \(a=48\lceil\log_2m\rceil\ge48\).
Jensen's inequality bounds the expected root-mean-square distance
by the left-hand side of \eqref{eq:arbitrary-rms}.
The expected-distance lower bound for \(T\ge m^2\), from arbitrary
initial supply locations and even with perfect foresight, is due to
Kanoria
\cite[Theorem~3 and its proof]{Kanoria}.
The replacement phase can maintain the tree counts and supply lists
as supplies arrive and depart, so the \(O(\log m)\) operations per
period and \(O(m)\) memory bounds from Section~\ref{sec:policy} still hold.

\subsection{Guarantees for all $T$ under balanced initial inventories}

\if
\begin{remark}[Stationary and limiting bounds]
\label{rem:other-initial}
For the harmonic-feedback policy in stationarity,
\eqref{eq:stationary-energy} from Proposition~\ref{prop:energy} and
\eqref{eq:transport-cost} from Lemma~\ref{lem:transport-energy} give
\begin{equation}\label{eq:stationary-cost}
 \sqrt{\E r_t^2}
 \le\frac1m+\sqrt{\frac53}\frac am
 <\frac{3a}{2m}\, .
\end{equation}
By Lemma~\ref{lem:symmetry}, the leaf-count chain converges to
stationarity from any initial inventory, so
\[
 \limsup_{t\to\infty}\sqrt{\E r_t^2}\le\frac{3a}{2m}\, ,
\] 
holds for arbitrary initial inventory. 
\end{remark}
\fi

\begin{corollary}[Balanced initial inventory]\label{cor:balanced}
For every \(m\ge2\), there is an initial inventory law for which the
harmonic-feedback policy, run without the initial replacement phase,
satisfies
\[
\biggl(\frac1T\sum_{t=1}^T\E r_t^2\biggr)^{1/2}
\le\frac1m+\sqrt{\frac53}\frac am
\qquad\text{for every }T\ge1\, .
\]
When \(n\) divides \(m\), equal counts in all leaves suffice.
\end{corollary}
\begin{proof}
Split the \(m\) supplies recursively as evenly as possible, using
independent fair coins to choose the larger child whenever a count is odd.
At each level the counts differ by at most one. Concavity of \(F_a\)
implies that this configuration maximizes
the contribution of every level, hence maximizes \(\Phi\).
Its count law is tree invariant. Thus \(\Phi_T-\Phi_0\le0\), and
\eqref{eq:master-energy} from the proof of Proposition~\ref{prop:energy}
bounds $\frac1T\sum_{t=1}^T\E G(t)$ %
by
\(20a^2/m^2\) for every \(T\). 
Choose \(J\) uniformly from \(\{1,\ldots,T\}\), independently
of the process. The count law at period \(J\) is a mixture
of the tree-invariant count laws at the individual periods,
and is therefore tree invariant. Applying
\eqref{eq:transport-cost} from Lemma~\ref{lem:transport-energy}
to this mixture gives
\[
\biggl(\frac1T\sum_{t=1}^T\E r_t^2\biggr)^{1/2}
\le \frac1n+\sqrt{\frac1{12T}\sum_{t=1}^T\E G(t)}
\le \frac1m+\sqrt{\frac53}\frac am\, .
\]
\end{proof}

\section{Design principles}\label{sec:heuristic}

In this section, we motivate key choices in the design of our policy and analysis: the potential function, the feedback rule, the correction $B(h,t)$ for inherited drift in potential, and the choice of $a$. 

\paragraph{Choice of potential function $\Phi$.}
As discussed in Section~\ref{sec:intro}, the analogy with two-dimensional static matching suggests that different length scales contribute comparably to matching cost in our one-dimensional dynamic model.
Accordingly, we would like to score the distribution of supply inventory in a way that treats different length
scales alike. This motivates our use of a logarithmic potential; a given percentage drop in supply
produces the same change in log inventory in a short interval as in a
long one.

We use the regularized discrete version of the logarithm
\(F_a(N)=\sum_{j=1}^{N}\frac{1}{j+a/2}\), with \(F_a(0)=0\). Removing one of \(N\ge1\) units
decreases \(F_a(N)\) by exactly the harmonic term \(1/(N+a/2)\). For an interval of length
\(p\), the corresponding quantity \(h=p/(N+a/2)\) approximates the
reciprocal of the supply density when \(N\gg a\), and remains finite when
the interval is empty.

What weight should the potential assign to an interval of length \(p\)?
If the interval contains no supply, demands arriving inside it contribute
at least \(\int_0^p\min\{z,p-z\}\,dz=p^2/4\) to the expected matching
distance per period, even with supplies immediately outside both
endpoints. This suggests weighting an interval by its squared length, leading to
our definition of the potential $\Phi$ in~\eqref{eq:Phi-D}.

At a well-stocked level, an interval of length \(p\) holds about
\(mp\gg a\) units, so the curvature of \(F_a\) is of order \((mp)^{-2}\).
Multiplying by the interval weight \(p^2\) and the supply arrival
probability \(p\) gives a curvature-loss bound of order
\(p^2\cdot p\cdot(mp)^{-2}=p/m^2\) per interval in~\eqref{eq:curvature-error}.
Summing over the \(1/p\) intervals gives a bound of order \(1/m^2\),
independent of the length scale.

\paragraph{Choice of feedback rule $b$.}
Consider an interval of length \(p\) whose children both contain supply.
Write \(h_+\) and \(h_-\) for the inverse inventories of the fuller and
less full children, respectively. Each child has weight \((p/2)^2\), so
deleting one unit reduces its weighted potential by \((p/2)h_+\) in the
fuller child and by \((p/2)h_-\) in the less full child. Relative to an
even split, a bias \(b\) shifts deletion probability \(b/2\) toward the
fuller child.
The resulting gain in expected potential is therefore
\(pb(h_--h_+)/4\).

The gain in potential is linear in \(b\).
Recall from Lemma~\ref{lem:transport-energy} that under tree symmetry,
the \(pb^2/12\) terms sum in expectation to the expected squared distance
from demand to proxy supply.
Charging the same price for
squared displacement at every level, and writing that price as
\(3/(2a)\), gives the local objective
\[
  \frac{pb}{4}(h_--h_+)-\frac{pb^2}{8a}.
\]
Its maximizer is \(b=a(h_--h_+)\), the feedback rule in~\eqref{eq:policy} for cases where both children are occupied.
When one child is empty, the occupied child must supply the match.

\paragraph{Correcting for inherited drift in potential: the choice of $B$.}
In choosing \(b\), we held the parent's deletion probability \(q\) fixed
and optimized how to divide it between the children. The parent's ancestors
determine \(q\). When the children have equal inventories,
\(b=0\), but if \(q>p\), both children lose inventory in expectation.
Their contribution to the linearized potential drift $D$ defined in \eqref{eq:Phi-D} 
is then negative.

Write \(t=(p-q)/a\), and use primes for quantities at the children. The
contribution of a split to \(D/a\) is \(p\,\E_{\mathrm{ch}}[t'h']\),
where \(\E_{\mathrm{ch}}\) averages equally over the two children. We
want the sum of these contributions to control
\(G/a^2=\sum_{v\in\T\setminus\T_{\ell_0}}p_v(b_v/a)^2\), where \(G\)
is the transport energy defined in~\eqref{eq:H-G}.
The equal-inventory example shows why we cannot ask each raw child
contribution to do this on its own. We seek a function \(B(h,t)\) for
which
\[
  \E_{\mathrm{ch}}[t'h']+\E_{\mathrm{ch}}[B(h',t')]-B(h,t)
  \;\ge\; \kappa (b/a)^2
\]
at every split, for some constant \(\kappa>0\). The difference form is
what makes \(B\) useful. After multiplying by \(p\) and summing over the
tree, each interior node's correction appears once as a child and once
as a parent, with opposite signs. Only the root and leaf corrections
survive.

How might one choose \(B\)? With counts and \(a\) fixed, a common
rescaling of \(p\) and \(q\) rescales \(h,t,h',t'\), and \(b/a\) by the
same factor. Both \(t'h'\) and \((b/a)^2\) scale as the square of that factor. We
therefore look for \(B\) among homogeneous quadratic polynomials in
\(h\) and \(t\).

To see which quadratic to try, consider two occupied children. Put
\(H=(h_-+h_+)/2\) and \(w=b/(2a)\). The feedback rule gives the two
child pairs
\[
  (h',t')=(H+w,t/2+w)
  \qquad\text{and}\qquad
  (h',t')=(H-w,t/2-w).
\]
Consequently, \(\E_{\mathrm{ch}}[t'h']=Ht/2+w^2\). The \(w^2\) term is
the gain we want to retain. In each child, \(h'\) and \(t'\) contain the
same signed \(w\) term, so \(h'-t'=H-t/2\) in both children. If we give \(B(h,t)\)
a factor \(h-t\), its other factor is linear. In averaging the child
corrections, the first factor is constant and the \(+w\) and \(-w\) terms
in the linear factor cancel. The correction thus leaves the favorable
\(w^2\) term intact.

For the second factor, consider an empty child. There \(t'=h'/2\), so
its raw contribution is already positive, \(t'h'=(h')^2/2\). A natural
choice is to set the correction to zero there, retaining that contribution
in full. The linear factor that vanishes when \(t=h/2\) is proportional
to \(h-2t\). These two choices give the candidate
\[
  B(h,t)=-c(h-t)(h-2t).
\]
For \(c>0\), the feasibility invariant \(t\le h/2\) from
Lemma~\ref{lem:invariants} also gives the desired sign \(B\le0\) at the
leaves.

To choose \(c\), return to the equally stocked split that exposed the
problem. When the parent is well stocked, each child has approximately
the parent's inverse inventory, \(h'\approx h\), and inherits \(t'=t/2\).
Setting \(h'=h\) in this calculation gives
\[
  \frac{th}{2}+B(h,t/2)-B(h,t)
  =\frac{1-3c}{2}\,th+\frac{3c}{2}\,t^2.
\]
Choosing \(c=1/3\) cancels the term linear in \(t\) and leaves \(t^2/2\),
which is nonnegative for either sign of the inherited drift. This motivates
the definition of $B(h,t)$ in~\eqref{eq:Bellman}. 
Lemma~\ref{lem:local} proves the local inequality with \(\kappa=1/30\)
for every possible configuration of parent and child inventories,
including splits where only one child contains supply.

\paragraph{The choice of \(a\) and the contribution of each scale.}
The drift estimates explain why we choose the feedback strength to be
proportional to the tree depth \(\ell_0\). After summing the local inequalities, the lower bound on
linearized drift in~\eqref{eq:aggregate}, from the proof of
Proposition~\ref{prop:energy}, has coefficient \(1/(30a)\) on the
transport energy \(G\). The upper bound on curvature loss
in~\eqref{eq:R-bound}, obtained by summing the second inequality of
Lemma~\ref{lem:local} across levels, has coefficient \(3\ell_0/(4a^2)\)
on that same energy. Taking \(a=48\ell_0\) leaves a positive coefficient
on \(G\) in the expected drift of the harmonic potential.

We now give a heuristic estimate of how much one well-stocked level
contributes to \(\E G\). At a node of length \(p\), use spatial signs for
the count imbalance \(\Delta=N_{\mathrm L}-N_{\mathrm R}\) and the bias
\(b=q_{\mathrm L}-q_{\mathrm R}\). In steady state, the node’s expected inventory change is zero, so \(\E q=p\) holds.
The conditional drift of \(\Delta^2\) is \(-2\Delta b+p+q\), so
we further obtain the exact identity \(\E[\Delta b]=p\) in steady state.

Now consider a scale with \(mp\gg a\), and suppose the node is typically
well stocked and nearly balanced, with \(N\approx mp\) and
\(|\Delta|\ll N\). The feedback formula then gives
\(b\approx 2a\Delta/(m^2p)\). Using this approximation in the steady state
identity predicts
\[
  \E[b^2]\approx\frac{2a}{m^2p}\,\E[\Delta b]=\frac{2a}{m^2}\, .
\]
Since the node lengths at a level sum to one, we estimate that each level should
contribute about \(2a/m^2\) to \(\E G\), or \(a/(6m^2)\) to expected squared
proxy displacement. With \(a\) of order \(\log m\) and order \(\log m\)
well-stocked levels, this gives the squared-distance scale
\((\log m)^2/m^2\). Proposition~\ref{prop:energy} bounds the total expected
transport energy over the entire tree, including the finest levels where
children can be empty.

\section{Discussion}\label{sec:discussion}\label{sec:taxi}

Our result can be interpreted as an average case result for the hard $k$-taxi 
problem, where only empty travel is charged. 
With \(k\ge2\) taxis,
mutually independent uniform pickups and destinations on \([0,1]\),
and travel times ignored, each ride has the matching dynamics with
\(m=k\): it replaces the dispatched
taxi's position by its destination. Thus, from arbitrary initial taxi
positions and for \(T\ge k^2\), Section~\ref{sec:finite} gives expected
total empty travel within a constant factor of the expected clairvoyant
optimum. Also note that proactive empty travel cannot improve that optimum: preserving
each taxi's assignments, the triangle inequality allows it to wait at
its initial position or previous dropoff, then travel directly to its
next pickup.

\paragraph{Use of AI.}
The author used GPT 6 Astra in developing the proof, through a
 back-and-forth over two months, and made a seminar presentation along the way which helped to clarify and refine the thinking, ultimately leading to the policy and proof presented here. 
Providing Astra with guidance to draw inspiration from Shor's conceptually related work~\cite{Shor1991} seemed to help it reach the new result here, 
in addition to the extensive background and intuition already provided in \cite{Kanoria}. 
Astra's first proof required computer verification with polynomial sign checks
of more than 600 coefficients of a polynomial. The author repeatedly parsed the proof returned and suggested directions to simplify and improve it, including joint design of the policy and analysis. The second proof
reduced the largest such check to just over 100 coefficients. The third
proof could be checked by hand but was still lengthy and somewhat unnatural. The fourth
further shortened the boundary calculations. The fifth proof introduced a particularly natural form of the
policy and
analysis, presented here. GPT 6 Astra was also used to create the Lean formalization, and to polish the writing.
In a separate trial, Fable~5, given only a pointer to the 
conjecture presented in 
\cite{Kanoria}, with no additional
guidance, eventually reported that the conjecture was likely false after examining
the performance of two policies it constructed. Earlier AI models from both Open AI and Anthropic made
no progress on the conjecture in the author's attempts, despite
substantial guidance.

Overall, the author was pleased with a particularly elegant and intuitive resolution to a fairly hard open problem, leveraging a Lyapunov drift approach which may not be typical in the context of spatial matching, motivating him to write it up and share it. 

\bibliographystyle{amsplain}
\bibliography{references}
\end{document}